\documentclass[12pt,letterpaper,reqno]{amsart}

\usepackage{amsmath}
\usepackage{amsfonts}
\usepackage{amssymb}
\usepackage{amsthm}
\usepackage{amsrefs}
\usepackage{amscd}
\usepackage{enumitem}
\usepackage[colorlinks=true]{hyperref}
\usepackage[nosort,noabbrev,capitalise]{cleveref}

\usepackage{latexsym}
\usepackage{mathrsfs}
\usepackage{psfrag}
\usepackage[dvips]{epsfig}
\usepackage{epsfig}
\usepackage[all]{xy}         

\newcommand{\mcC}{\mathcal{C}}

\newcommand{\mcP}{\mathcal{P}}

\swapnumbers 

\theoremstyle{plain}                              
\newtheorem{thm}{Theorem}[section]
\newtheorem{prop}[thm]{Proposition}
\newtheorem{defn}[thm]{Definition}
\newtheorem{lem}[thm]{Lemma}

\theoremstyle{definition}                         

\theoremstyle{remark}                             

\numberwithin{equation}{section}

\newcommand{\R}{\mathbb{R}}                     

\begin{document}

\title[Space-Filling Curves]{Space-Filling Curves}

\author{Shihan Kanungo}

\address{Department of Mathematics, San Jos\'e State University, San
  Jos\'e, CA 95192-0103}

\email{shihankanungo@sjsu.edu}

\thanks{\textsc{Math 231A: Real Analysis I}}

\date{\today}

\begin{abstract}
    We examine space-filling curves, which are surjective continuous maps from $[0,1]$ to some higher-dimensional space, usually the unit square $[0,1]^2$. In particular, we define Peano's curve and Lebesgue's curve, and state some of their properties. We also discuss the Hahn-Mazurkiewicz theorem, which characterizes those subsets of $\R^n$ that are the image of a space-filling curve. Finally, we discuss real-world applications of Hilbert curves, in particular Google's $S2$ Cells.
\end{abstract}

\maketitle

\section{Introduction}  \label{sec:intro}
A space-filling curve is a surjective continuous map $f: [0,1]\to [0,1]^n$ for some $n\ge 2$, which can be thought of as a curve inside a unit square that somehow manages to hit every single point. Since a curve has no thickness, it is extremely counterintuitive that it could fill a square. The first indication that this could even be possible came in the mid 19th century, when Georg Cantor (1845 -- 1918) proved that $[0,1]$ and $[0,1]^2$ have the same cardinality, i.e. there is a way to bijectively map the unit line segment to the unit square. However, his map flings points that are very close in $[0,1]$ to far-away points in $[0,1]^2$; so is not continuous. Even so, this result was already so unbelievable that in his 1877 correspondence of this result to his friend and fellow mathematician, Richard Dedekind, Cantor remarked, ``I see it, but I don’t believe it!''. Naturally, the next question to ask was: is there a \textit{continuous} surjective map $[0,1]\to [0,1]^2$. In other words, is there a space-filling curve? 

\begin{figure}[h!]
    \centering
    \includegraphics[width=0.45\linewidth]{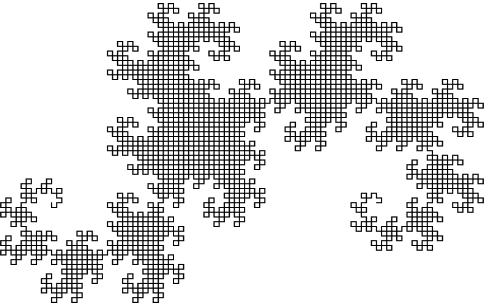}
    \caption{The Dragon Fractal}
    \label{fig: dragon fractal}
\end{figure}

Motivated by Cantor's work, Giuseppe Peano (1858 -- 1932), an Italian mathematician and professor at the University of Turin, first discovered a space-filling curve in 1890 by constructing a surjective continuous function from $[0,1]$ to the unit square $[0,1]^2$. Soon after, many other space-filling curves were discovered, including the Hilbert curve (a variant of Peano's curve), Lebesgue's curve, and the famous example of the \textit{dragon fractal}, shown in \cref{fig: dragon fractal}.

In this paper, we will look at Hilbert's simplified version of the Peano curve, and Lebesgue's curve, another space filling curve $[0,1]\to [0,1]^2$. We also discuss the Hahn-Mazurkiewicz Theorem, which classifies which subsets of $\R^n$ can be the image of some continuous map $f:[0,1]\to\R^n$.

\section{Basic concepts}\label{sec: basic concepts}

In this section, we define the two important space-filling curves we will be studying: Peano's curve and Lebesgue's curve.

\subsection{Peano's curve}
\begin{figure}[h!]
    \centering
    \includegraphics[width=0.75\linewidth]{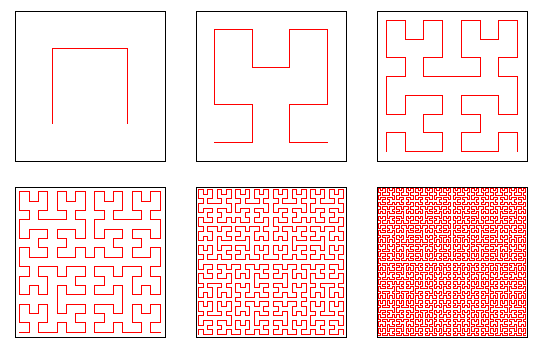}
    \caption{The Hilbert Curve}
    \label{fig: Hilbert curve}
\end{figure}

There are many variants of Peano's curve. The one we present is the simplest variant, called the \textit{Hilbert curve}. We will give  a visual description of the Hilbert curve, which is constructed iteratively as shown in \cref{fig: Hilbert curve}.

The idea of the construction is similar to the construction of the Cantor set. We let \textit{quartic intervals of generation $k$} be the intervals obtained when $[0,1]$ is partitioned into $4^k$ intervals of equal length. We let $I^k$ denote the set of quartic intervals of generation $k$. So there are $4^k$ intervals, each of length $4^{-k}$, in $I^k$, and they all have the form $[\tfrac{i}{4^k},\tfrac{i+1}{4^k}]$.

Each interval in $I^k$ is the union of four intervals in $I^{k+1}$. Conversely, every interval in $I^k$ is a subset of a unique interval in $I^{k-1}$. Generalizing this, we make the following definition.
\begin{defn}\label{defn: chain of quartic intervals}
    A \emph{chain} of quartic intervals is a sequence $(I_i)_{i=0}^\infty$ of intervals with $I_i\in I^i$ such that
    \[I_0\supset I_1\supset I_2\supset\cdots.\]
\end{defn}
We have the following elementary properties related to chains.
\begin{prop}\label{prop: chains of quartic intervals} Chains of quartic intervals satisfy the following properties:
\begin{enumerate}[label = $(\alph*)$]
    \item If $(I_i)^\infty_{i=0}$ is a chain, then there exists a unique $t\in [0,1]$ such that $t\in I_i$ for all $i$,
    \item For each $t\in [0,1]$ there is a chain $(I_i)^\infty_{i=0}$ satisfying $(a)$, and
    \item For all but countably many $t$, the chain in $(b)$ is unique.
\end{enumerate}
\end{prop}
Chains are related to base-4 representation in the following way: Each quartic interval $I\in I^k$ consists of all numbers $x$ that have a base-4 expansion starting with $\big(\overline{0. a_1a_2\dots a_k}\big)_4$, where $a_1,\ldots,a_k$ are fixed digits between $0$ and $3$. So we can think of each $I_k$ in the chain as fixing the $k$-th digit in the base-4 expansion of $t$. Notice that this argument is similar to one used to prove that the Cantor sets consists of those numbers with no $1$ in their base-3 expansion.

Now, the Hilbert curve will be constructed iteratively, with each iteration mapping intervals in $I^k$ to a $2$-dimensional analog of quaternary intervals.

The \textit{dyadic squares} of generation $k$ are the set of $4^k$ squares obtained when $[0,1]^2$ is partitioned into a $2^k\times 2^k$ grid of $2^{-k}\times 2^{-k}$ squares. We let $S^k$ denote the set of dyadic squares of generation $k$. Mirroring \cref{defn: chain of quartic intervals}, we define
\begin{defn}
    A \emph{chain} of dyadic squares is a sequence $(S_i)_{i=0}^\infty$ with $S_i\in S^i$ such that
    \[S_0\supset S_1\supset S_2\supset\cdots.\]
\end{defn}
We pause here to mention that we are assuming all dyadic squares, as well as all quartic intervals to be closed.

There is an analog of \cref{prop: chains of quartic intervals} for dyadic squares.
\begin{prop}\label{prop: chains of dyadic squares} Chains of dyadic squares satisfy the following properties:
\begin{enumerate}[label = $(\alph*)$]
    \item If $(S_i)^\infty_{i=0}$ is a chain, then there exists a unique $x\in [0,1]^2$ such that $x\in S_i$ for all $i$,
    \item For each $x\in [0,1]^2$ there is a chain $(S_i)^\infty_{i=0}$ satisfying $(a)$, and
    \item For all but countably many $x$, the chain in $(b)$ is unique.
\end{enumerate}
\end{prop}

Now, we have the following important definition.
\begin{defn}\label{defn: dyadic correspondence}
    A \emph{dyadic correspondence} is a bijective map $\phi$ from the set of quartic intervals to the set of dyadic squares such that $\phi(I^k)=S^k$ and if $I\supset J$ then $\phi(I)\supset \phi(J)$.
\end{defn}
Note that these exist because $|I^k| = |S^k| = 4^k$. Now, from each dyadic correspondence, we can get a map from $[0,1]$ to $[0,1]^2$ as follows: for each $t\in [0,1]$ there is a chain $(I_i)_{i=0}^\infty$ with $t\in I_i$ for all $i$. Since $\phi$ respects inclusion, $(\phi(I_i))_{i=0}^\infty$ is also a chain of dyadic squares. So there is a unique $x\in [0,1]^2$ with $x\in \phi(I_i)$ for all $i$. We define $\mcP(t)=x$. Since $\phi$ is a bijection, it follows that $\mcP([0,1])$ differs from $[0,1]^2$ by a countable set (as countably many $t\in [0,1]$ have multiple chains corresponding to $t$).

However, we don't know whether $\mcP(t)$ is continuous or not. We will show in the following section that there exists a dyadic correspondence such that the corresponding curve $\mcP$ is continuous, and this is the Hilbert curve.
\subsection{Lebesgue's curve}
Lebesgue's space-filling curve is defined using the Cantor set $\mcC$, which contains all elements of $[0,1]$ which don't have a $1$ in their ternary expansion. Then, in base-3, each element of $\mcC$ can be written as $\overline{0.(2a_1)(2a_2)(2a_3)\dots}_3$ where $a_i\in \{0,1\}$. Recall that the Cantor-Lebesgue function $L: \mcC\to [0,1]$ is defined as
\begin{align}\label{eqn: Cantor-Lebesgue}
L(\overline{0.(2a_1)(2a_2)(2a_3)\dots}_3) = \overline{0.a_1a_2a_3\dots}_2
\end{align}
for elements of $\mcC$ and then extending linearly to all of $[0,1]$. Then $L$ is an example of a continuous function which has derivative $0$ almost everywhere, but goes from $0$ to $1$.

\begin{figure}
    \centering
    \includegraphics[width=0.5\linewidth]{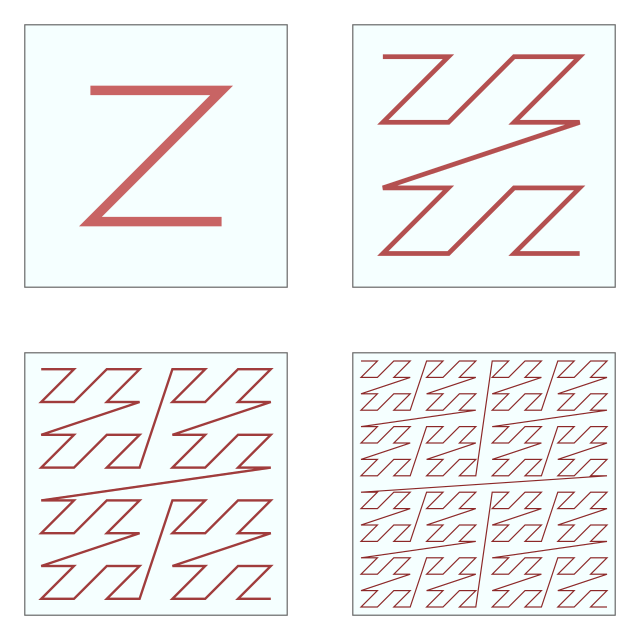}
    \caption{The Lebesgue Curve}
    \label{fig: lebesgue curve}
\end{figure}

Similarly, we can define $\ell: \mcC \to [0,1]^2$
\begin{align}\label{eqn: Lebesgue curvee}
\ell(\overline{0.(2a_1)(2a_2)(2a_3)\dots}_3) = \left(\overline{0.a_1a_3a_5\dots}_2,\overline{0.a_2a_4a_6\dots}_2\right),
\end{align}
which is clearly surjective. Morever, similarly to the proof that $L$ is continuous, we have:
\begin{prop}
    The function $\ell: \mcC\to [0,1]^2$ is continuous.
\end{prop}
Finally, we can extend $\ell$ linearly to a function $\widetilde{\ell}: [0,1]\to [0,1]^2$. More formally, for each gap $(a_n,b_n)$ in the Cantor set $\mcC$, we define
\[\widetilde{\ell}(x) = \ell(a_n) + \frac{x_n-a_n}{b_n-a_n}(\ell(b_n)-\ell(a_n))\]
for $x\in [a_n,b_n]$. Just like the extension of $L$ to $[0,1]$, this linear extension is also continuous.
\begin{prop}
    The function $\widetilde{\ell}:[0,1]\to[0,1]^2$ is continuous.
\end{prop}
This function $\widetilde{\ell}$ is Lebesgue's space-filling curve, which can also be described iteratively like the Hilbert curve. This is shown in \cref{fig: lebesgue curve}.

\section{Main results}\label{sec: main results}
In this section, we look at the existence and some properties of the Hilbert Curve. Then, we state the Hahn-Mazurkiewicz theorem and discuss the ideas of the proof.
\subsection{Properties of the Hilbert curve}
It turns out that there exists a choice of the dyadic correspondence $\phi$ such that the corresponding map $\mcP: [0,1]\to [0,1]^2$ satisfies the following nice properties.
\begin{thm}\label{thm: Peano curve}
    There exists a continuous and surjective map $\mcP: [0,1]\to[0,1]^2$ such that the following two properties holds
    \begin{enumerate}[label = $(\alph*)$]
        \item $\mcP$ satisfies a Lipschitz condition of exponent $1/2$:
        \begin{align}\label{eqn: lipschitz condition}
        \lVert\mcP(t)-\mcP(s)\rVert\le M|t-s|^{1/2}
        \end{align}
        for some constant $M$.
        \item For any $[a,b]\subset [0,1]$, $\mcP([a,b])$ is a compact subset of $[0,1]^2$ that has Lebesgue measure of exactly $b-a$.
    \end{enumerate} 
\end{thm}
Then, we can remove sets of measure zero from $[0,1]$ and $[0,1]^2$ such that $\mcP$ becomes a bijection, and then we have that $E\subseteq [0,1]$ is measurable if and only if $\mcP(E)\subseteq [0,1]^2$ is. Additionally, $m(E)=m(\mcP(E))$. It follows that from a measure-theoretic point of view, $[0,1]$ and $[0,1]^2$ are ``isomorphic'' in some sense. We can also use space-filling curves like these to turn an integral over $[0,1]^2$ to a one-dimensional integral over $[0,1]$. In practice however, evaluating such a integral would not be easy.

\begin{figure}[h!]
    \centering
    \includegraphics[width=0.7\linewidth]{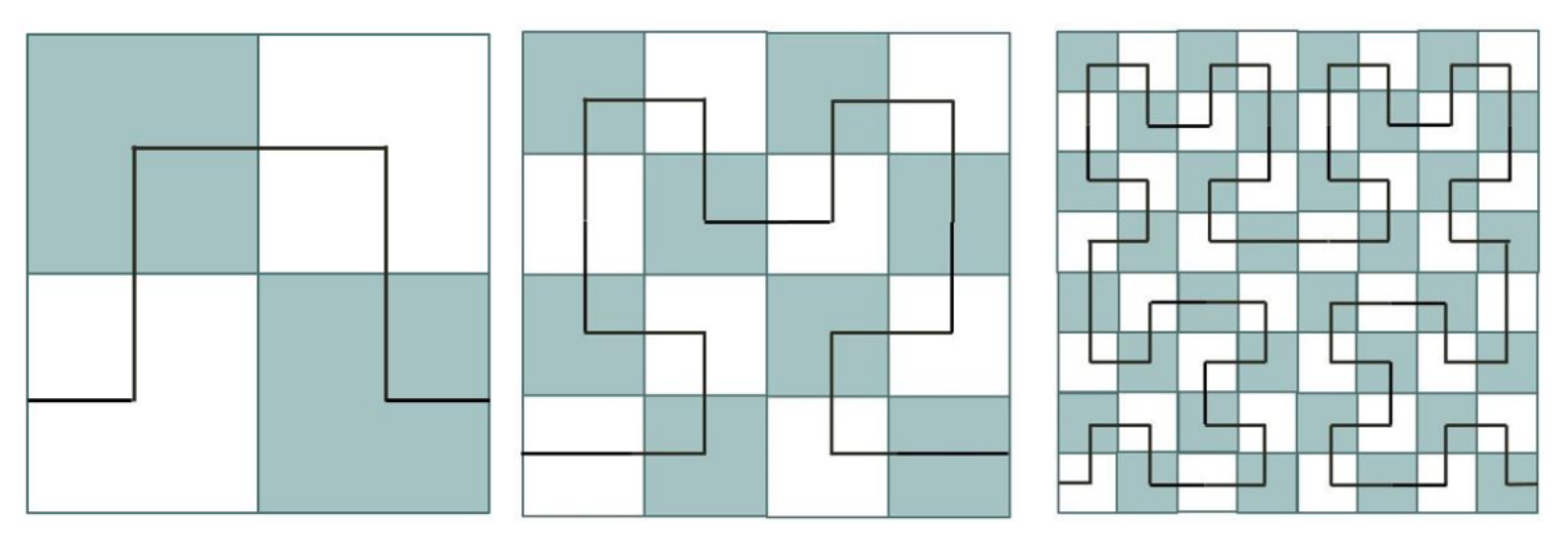}
    \caption{Construction of $\phi$}
    \label{fig: dyadic correspondence}
\end{figure}

To prove \cref{thm: Peano curve}, we use the following result.
\begin{lem}
    There is a unique dyadic correspondence $\phi$ such that
    \begin{enumerate}[label = $(\alph*)$]
        \item If $I$ and $J$ are two adjacent intervals in the same generation, then $\phi(I)$ and $\phi(J)$ are adjacent squares as well.
        \item If $I_{-}$ is the left-most interval in generation $k$ and $I_+$ is the right-most, then $\phi(I_{-})$ is the bottom-left-most square in generation $k$ and $\phi(I_+)$ is the bottom-right-most.
    \end{enumerate}
\end{lem}
The explicit construction for $\phi$ is shown in \cref{fig: dyadic correspondence}. 

To specify $\phi$, we simply need to give an ordering of the squares in $S^k$. Then the $i$'th interval in $I^k$ maps to the $i$'th square in the ordering. In the figure, we start at the bottom left and then follow the line to get the order of the squares till we end up at the end at the bottom right (Note how this is similar to \cref{fig: Hilbert curve}, but with the dyadic squares overlaid under the curve). 

\medskip 
With the square (and all its continuous images) revealed as a continuous image of $[0, 1]$, the question arose as to the general structure of a continuous image of a line segment. In 1908, A. Schoenflies found a criterion which, by its very nature, only applies to subsets of the two-dimensional plane and never entered the mathematical mainstream. The further pursuit of this question led to the new topological concept of local connectedness and to a complete answer found independently by S. Mazurkiewicz and H. Hahn in 1913. It turned out that a set is the continuous image of a line segment if and only if it is compact, connected, and locally connected, and this criterion applies not only to $\mathbb{R}^n$ but more generally to Hausdorff spaces.

\section{The Hahn-Mazurkiewicz Theorem}
Let's start by defining the key concepts related to the Hahn-Mazurkiewicz theorem.
\begin{defn}
    A subset $A$ of $\R^n$ is \emph{compact} if every sequence has a convergent subsequence with limit in $A$. By the Bolzano-Weirstrass theorem, this is equivalent to $A$ being closed and bounded. 
    We say that $A$ is \emph{connected} if there does not exist non-empty disjoint $U,V\in \R^d$ such that $A = A\cap (U\sqcup V)$. 
    Finally, we say $A$ is \emph{weakly locally-connected} if for all $a\in A$ and $\eta>0$, there exists $0<\varepsilon<\eta$ such that if $x\in B_{\varepsilon}(a)$, then $x$ and $a$ lie in the same connected component of $A\cap B_\eta(a)$.
\end{defn}
It turns out that weakly locally-connectedness is equivalent to the more famous notion of locally-connectedness. However, weakly locally-connectedness at a point is \textit{not} the same as locally connectedness at a point. Now, we may state the Hahn-Mazurkiewicz theorem.
\begin{thm}[Hahn-Mazurkiewicz]
    A subset $A\subseteq \R^n$ is the image of some continuous map $f: [0,1]\to \R^n$ if and only if $A$ is compact, connected, and weakly locally-connected.
\end{thm}
So this theorem classifies what kinds of sets a space-filling curve can occupy. Proving that the image of any space filling curve $f:[0,1]\to \R^n$ is compact, connected, and weakly locally-connected is relatively straightforward. To prove the other direction, we use ideas coming from the construction of Lebesgue's curve.

First, we have Hausdorff's theorem, which is stated as follows.
\begin{thm}[Hausdorf]
    Every compact subset of $\R^n$ is the continuous image of the Cantor set.
\end{thm}
To prove this, we create some sort of analog of the dyadic squares in $[0,1]^2$. More precisely, we have a series of finer and finer covers of our compact set $A$ such that each ball in some cover is a subset of a ball in a previous cover. Then, we can describe each element of $A$ by a chain of increasingly smaller balls in the covers. Then, for each base-3 representation of $x\in \mcC$ we can assign an element of $A$, which gives a surjective map, which is also continuous. This process is similar to \eqref{eqn: Cantor-Lebesgue} and \eqref{eqn: Lebesgue curvee}. 

So we have a continuous and surjective map $g:\mcC\to A$. Then the final step in the proof is to extend this linearly to $[0,1]$. The idea is that for each gap $(a,b)$ in the Cantor set, instead of simply extending our map linearly between $g(a)$ and $g(b)$ like we did for the Lebesgue curve, we instead find a (continuous) path from $g(a)$ to $g(b)$ that always stays inside our set $A$. The existence of such a path is made possible by the conditions that $A$ is connected and weakly locally-connected. Then this process yields a continuous function $f: [0,1]\to A$, proving the theorem.

However, while we now know many surjective and continuous functions $[0,1]\to [0,1]^2$, and in general from $[0,1]$ to any compact, connected, and weakly locally-connected subset of $\R^n$, none of them are injective. It turns out that this is for good reason---it's impossible! 
\begin{thm}[Netto]
    Any bijective map $f: [0,1]\to [0,1]^2$ is necessarily discontinuous.
\end{thm}

Of course, the same result holds if we replace $[0,1]^2$ with $\R^n$ or $[0,1]^n$ for $n\ge 2$.

\section{Applications} \label{sec:applications}
At first, space-filling curves might seem simply like a nice mathematical construction with some theoretical implications. It definitely doesn't seem like that they have any use in the real world, as after all, it's impossible to make a physical copy of one. But suprisingly, it turns out that they have an extremely important application: in Google maps!

In Google maps, the programmers want a way to have \textit{cache locality}. This means that there is a way to index all the positions on the earth such that if you move a little bit on the earth, then you only move a little bit in the memory. The reason that cache locality is so important is that it speeds up computing time immensely. In fact, it turns out that using the Hilbert curve, we can represent any square centimeter of the earth by a 64-bit integer! 

The formal name for this tool is ``$S2$ cells,'' which were developed by Eric Veach at Google. The idea is as follows. First, we put the earth inside a massive cube. Then, we put a Hilbert curve on each face of the cube such that the ends of curves on adjacent faces connect. 

\begin{figure}[h!]
    \centering
    \includegraphics[width=0.5\linewidth]{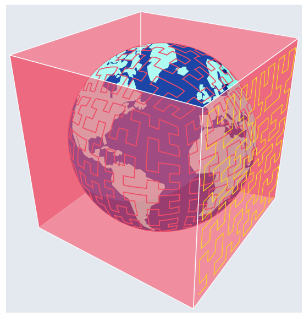}
    \caption{Google's $S2$ cells}
    \label{fig: s2 cells}
\end{figure}

Thus, we get a continuous function from a line segment to the surface of the cube, depicted in \cref{fig: s2 cells}. Next, we use this to create a continuous function from the line segment to the surface of the earth as follows: for a point $P$ on the surface of the earth, draw the ray from the center $0$ of the earth to $P$, and it will intersect the cube at a point $P'$. Then, this point $P'$ corresponds to some point on the line segment, and we assign this point to $P$.

Since the Hilbert curve is continuous, this new map is also continuous. Additionally, if two points are close together on the Hilbert curve, then their corresponding points are close together on the line segment---indeed, this follows from the Lipschitz condition \eqref{eqn: lipschitz condition}.

Using the dyadic squares of generation $k$, we can divide the surface of the earth into $6\cdot 4^k$ squares ($4^k$ squares from each of the $6$ squares), and it only takes $2k+4$ bits to encode this: $3$ bits to specify a number from $0$ to $5$, corresponding to one of the $6$ faces, and then $2k$ bits to encode one of the $4^k$ squares, and finally a ``1'' at the end to encode that there are no more bits.

Now, we can use the Hilbert curve (and the corresponding dyadic correspondence) to pick the order of the squares so that close-by squares correspond to close-by numbers in the encoding. To do this, note that each quartic interval is $I=[\tfrac{i}{4^k},\frac{i+1}{4^k}]$, where in base-2, $i$ is $\overline{a_1a_2\dots a_{2k}}_2$ for $a_{j}\in\{0,1\}$. So then, we set the square corresponding to the sequence of bits $a_1a_2\dots a_{2k}$ to be $\phi(I)$, where $\phi$ is our dyadic correspondence shown in \cref{fig: dyadic correspondence}. 

Finally, from each dyadic square on the face, we can assign a region on the earth by projecting it down: for each $P'$ on the square, we send it to the point $P$ on the sphere such that $O,P,P'$ are collinear.

If we set $k=30$, then each cell is represented by a $2k+4=64$-bit integer. If we work out the computations, it turns out that each cell has area between $0.48\mathrm{cm}^2$ and $0.93\mathrm{cm}^2$, so we can represent any point on the earth with a $64$-bit integer, with accuracy better than a centimeter!

A similar application is how GPUs use $Z$-order curves (higher-dimensional generalizations of the Lebesgue curve) to code multi-dimensional data into one dimension, while preserving locality (points close together in high dimensions are also close together on the line). Just like $S2$ cells, this makes computation much easier.

Space-filling curves are also used in load balancing: if computational domains need to be more precise in certain locations, space filling curves can deal with this by refining the curve in those places. For example, we might only need the generation $8$ Hilbert curve for the top left, bottom left, and bottom right quadrants of the unit square, but in the top right we might need generation $15$ instead. For an example of this used in practice, see \cite{Liu2017LoadBU}.

\end{document}